\documentclass[reqno]{amsart}
\usepackage{amssymb}

\usepackage{amscd}
\usepackage{amsmath}
\usepackage{graphicx}

\newtheorem{theorem}{Theorem}[section]
\theoremstyle{plain}

\newtheorem{corollary}[theorem]{Corollary}

\newtheorem{definition}[theorem]{Definition}

\newtheorem{lemma}[theorem]{Lemma}

\newtheorem{proposition}[theorem]{Proposition}
\newtheorem{remark}[theorem]{Remark}

\renewcommand\Bbb{\mathbb}

\begin{document}
\title[Complex interpolation of compactness for lattice couples]{Complex interpolation of compact operators mapping into lattice couples.}
\author{Michael Cwikel}
\address{Department of Mathematics, Technion - Israel Institute of Technology, Haifa
32000, Israel}
\email{mcwikel@math.technion.ac.il}
\thanks{The research was supported by the Technion V.P.R.\ Fund and by the Fund for
Promotion of Research at the Technion.}
\maketitle

\begin{abstract}
Suppose that $(A_0,A_1)$ and $(B_0,B_1)$ are Banach couples, and
that $T$ is a linear operator which maps $A_0$ compactly into $B_0$
and $A_1$ boundedly (or even compactly) into $B_1$.

\noindent
Does this imply that $T:[A_0,A_1]_\theta \to [B_0,B_1]_\theta$ compactly?

\noindent
(Here, as usual, $[A_0,A_1]_\theta$ denotes the complex interpolation space of Alberto
Calder\'on.)

\noindent
This question has been open for 44 years. Affirmative answers are known
for it in many special cases.

\noindent
We answer it affirmatively in the case where $(A_0,A_1)$ is arbitrary and
$(B_0,B_1)$ is a couple of Banach lattices having absolutely continuous norms
or the Fatou property.
\end{abstract}

\smallskip \urladdr{http://www.math.technion.ac.il/\~{}mcwikel/}


\section{\label{intro}Introduction.}

All Banach spaces in this paper will be over the complex field. The closed
unit ball of a Banach space $A$ will be denoted by $\mathcal{B}_{A}$. For
any two Banach spaces $A$ and $B$, the notation $T:A\overset{b}{\rightarrow }%
B$ will mean, just like the usual notation $T:A\rightarrow B$, that $T$ is a
linear operator $T$ defined on $A$ (and also possibly defined on a larger
space) and it maps $A$ into $B$ boundedly. The notation $T:A\overset{c}{%
\rightarrow }B$ will mean that $T:A\overset{b}{\rightarrow }B$ with the
additional condition that $T$ maps $A$ into $B$ compactly.

\smallskip We will write $A\overset{1}{\subset }B$ when $A$ is continuously
embedded with norm $1$ into $B$, and $A\overset{1}{=}B$ when $A$ and $B$
coincide with equality of norms.

For each \textbf{\textit{Banach couple}} (or \textbf{\textit{interpolation
pair}}) $\vec{A}=(A_{0},A_{1})$ and each $\theta \in [0,1]$, we will let $%
[A_{0},A_{1}]_{\theta }$ denote the complex interpolation space of Alberto
Calder\'{o}n \cite{ca}. We also let $A_{j}^{\circ }$ denote the closure of $%
A_{0}\cap A_{1}$ in $A_{j}$ for $j=0,1$. The couple $(A_{0},A_{1})$ is
called \textbf{\textit{regular}} if $A_{j}^{\circ }=A_{j}$ for $j=0,1$. The
spaces $A_{0}\cap A_{1}$ and $A_{0}+A_{1}$ are Banach spaces when they are
equipped with their usual norms (as e.g., on p.\ 114 of \cite{ca}).

For any two fixed Banach couples $\vec{A}=(A_{0},A_{1})$ and $\vec{B}%
=(B_{0},B_{1})$, the notation $T:\vec{A}\overset{c,b}{\rightarrow }\vec{B}$
will mean that the linear operator $T:A_{0}+A_{1}\rightarrow B_{0}+B_{1}$
satisfies $T:A_{0}\overset{c}{\rightarrow }B_{0}$ and $T:A_{1}\overset{b}{%
\rightarrow }B_{1}$.
The notation $\vec{A}\blacktriangleright \vec{B}$ will mean that every
linear operator $T:A_{0}+A_{1}\rightarrow B_{0}+B_{1}$ which satisfies $T:%
\vec{A}\overset{c,b}{\rightarrow }\vec{B}$ also satisfies $%
T:[A_{0},A_{1}]_{\theta }\overset{c}{\rightarrow }[B_{0},B_{1}]_{\theta }$
for every $\theta \in (0,1)$. The notation $(*.*)\blacktriangleright \vec{B}
$ for some fixed Banach couple $\vec{B}$ will mean that $\vec{A}%
\blacktriangleright \vec{B}$ for every Banach couple $\vec{A}$. Analogously,
the notation $\vec{A}\blacktriangleright (*.*)$ for some fixed Banach couple
$\vec{A}$ will mean that $\vec{A}\blacktriangleright \vec{B}$ for every
Banach couple $\vec{B}$.

\smallskip Some forty-four years ago, Calder\'{o}n proved \cite{ca} that $%
(*.*)\blacktriangleright \vec{B}$ for all Banach couples $\vec{B}$ which
satisfy a certain approximation condition. Since then it has been
established that $\vec{A}\blacktriangleright \vec{B}$ for a large variety of
other different choices of $\vec{A}$ and $\vec{B}$. (See, e.g., the twelve
papers and website referred to on p.\ 72 of \cite{cj}, and \cite{cj}
itself.) However we still do not know whether $\vec{A}\blacktriangleright
\vec{B}$ holds for \textit{all} choices of $\vec{A}$ and $\vec{B}$. Let us
slightly extend the notation just introduced here and rephrase this by
saying that we are still trying to discover whether or not $%
(*.*)\blacktriangleright (*.*)$.

\smallskip In this paper we shall add to the library of known examples of
couples $\vec{A}$ and $\vec{B}$ satisfying $\vec{A}\blacktriangleright \vec{B%
}$ in the context of spaces of measurable functions. We shall use the
terminology \textbf{\textit{lattice couple}} to mean a Banach couple $\vec{A}%
=(A_{0},A_{1})$ where both $A_{0}$ and $A_{1}$ are complexified Banach
lattices of measurable functions defined on the same $\sigma $-finite
measure space.

\smallskip Cobos, K\"{u}hn and Schonbek (\cite{cks} Theorem 3.2 p.\ 289) proved that $%
\vec{A}\blacktriangleright \vec{B}$ whenever \textit{both} $\vec{A}$ and $%
\vec{B}$ are lattice couples, provided that $B_{0}$ and $B_{1}$ both have
the Fatou property, or that at least one of $B_{0}$ and $B_{1}$ has
absolutely continuous norm. Subsequently, Cwikel and Kalton (\cite{CwK}
Corollary 7 part (c) on p.\ 270) generalized this result by showing that $%
\vec{A}\blacktriangleright (*.*)$ for any lattice couple $\vec{A}$.

\smallskip In this paper we shall obtain a different generalization of the above
mentioned result of \cite{cks}, namely{\tiny \ }we will show that $%
(*.*)\blacktriangleright \vec{B}$ for every lattice couple $\vec{B}=(B_0,B_1)$
satisfying one or the other of the same conditions imposed by Cobos,
K\"{u}hn and Schonbek. In fact some other weaker conditions on $\vec{B}$ are
also sufficient. Roughly speaking, as indeed the reader might naturally
guess, our approach is to take the ``adjoint'' of the above mentioned result
$\vec{A}\blacktriangleright (*.*)$ of \cite{CwK}, using arguments in the
style of Schauder's classical theorem about adjoints of operators. But this
is apparently not quite as simple to do as one might at first expect. In
fact it will be convenient to use a somewhat more general ``abstract''
version of Schauder's theorem.

\smallskip As we will show in a forthcoming paper, our main result here is
one of the components needed to show that $(*.*)\blacktriangleright
\vec{B}$ also for certain other couples $\vec{B}$, \textit{non} lattice
couples which have not been considered previously in this context, and
which are in some sense quite close to the couple $\left( \ell ^{\infty
}(FL^{\infty }),\ell ^{\infty }(FL_{1}^{\infty })\right) $. This latter
couple is very important and interesting because (cf.\ Proposition 3 on p.\
356 of \cite{ckm}, or \cite{duke} pp.\ 339--340) determining whether or not
$(*.*)\blacktriangleright \left( \ell ^{\infty }(FL^{\infty }),\ell
^{\infty }(FL_{1}^{\infty })\right) $ is equivalent to determining whether
or not $(*.*)\blacktriangleright (*.*)$.

\smallskip Recently Evgeniy Pustylnik \cite{pustylnik} has obtained a very general
compactness theorem which in some ways is similar and in some ways more
general than ours. In fact his result applies to a wide range of
interpolation methods, not just Calder\'on's, and the spaces $B_0$ and
$B_1$ in his ``range" couple can be rather more general than Banach
lattices. On the other hand they have to satisfy certain conditions which
we do not require for our range spaces.

\pagebreak

\section{\protect\smallskip \label{rgaas}A rather general
Arzel\`{a}-Ascoli-Schauder theorem.}

\smallskip In this section we describe the result which will play the
r\^{o}le of Schauder's theorem for the proof of our main result.

Let us recall that a \textbf{\textit{semimetric space}} $\left( X,d\right) $%
, also often referred to as a \textbf{\textit{pseudometric space}}, is
defined exactly like a metric space, except that the condition $d(x,y)=0$
for a pair of points $x,y\in X$ does not imply that $x=y$. (However $%
d(x,x)=0 $ for all $x\in X$.) Each semimetric space $(X,d)$ gives rise to a
metric space $\left( \widetilde{X},\widetilde{d}\right) $ in an obvious way,
where $\widetilde{X}$ is the set of equivalence classes of $X$ defined by
the relation $x\sim y\Longleftrightarrow d(x,y)=0$.

Here are three definitions and three propositions concerning an arbitrary
semimetric space $(X,d)$. The definitions are exactly analogous to standard
definitions for metric spaces, and the propositions are proved exactly
analogously to the proofs of the corresponding standard propositions in the
case of metric spaces, or by invoking those standard propositions for the
particular metric space $\left( \widetilde{X},\widetilde{d}\right) $.

\begin{definition}
\label{zball}Let $B(x,r)$ denote the \textit{\textbf{ball of radius }}$r$
\textit{\textbf{centred at}} $x$, i.e., for each $x\in X$ and $r>0$, we set $%
B(x,r)=\{y\in X:d(x,y)\le r\}$.
\end{definition}

\begin{definition}
\label{tobo}The semimetric space $(X,d)$ is said to be \textbf{\textit{%
totally bounded}} if, for each $r>0$, there exists a finite set $%
F_{r}\subset X$ such that $X=\bigcup_{x\in F_{r}}B(x,r)$.
\end{definition}

\begin{definition}
\label{dsep}The semimetric space $(X,d)$ is said to be \textit{\textbf{%
separable}} if there exists a countable set $Y\subset X$ such that $%
\inf_{y\in Y}d(x,y)=0$ for each $x\in X$.
\end{definition}

\begin{proposition}
\label{tbsep}If $(X,d)$ is totally bounded, then it is separable.
\end{proposition}

\begin{proposition}
\label{phir} $(X,d)$ is not totally bounded if and only if for some $r>0$
there exists an infinite set $E\subset X$ such that $d(x,y)>r$ for all $%
x,y\in E$ with $x\ne y$.
\end{proposition}

\begin{proposition}
\label{koshy} $(X,d)$ is totally bounded if and only if every sequence $%
\left\{ x_{n}\right\} _{n\in \Bbb{N}}$ in $X$ has a \textit{\textbf{Cauchy
subsequence}}, i.e., a subsequence $\left\{ x_{n_{k}}\right\} _{k\in \Bbb{N}}
$ which satisfies $\lim_{N\rightarrow \infty }\sup \left\{
d(x_{n_{p}},x_{n_{q}}):p,q>N\right\} =0$.
\end{proposition}

\smallskip For the benefit of any reader who may happen to find them
helpful, we have allowed ourselves the luxury of including the (very
standard) proofs of Propositions \ref{tbsep}, \ref{phir} and \ref{koshy} in
an appendix (Section \ref{pps}). Since the writer and the reader(s) of this
paper can choose to keep it only electronic, we do not have to worry too
much about wasting paper.

The following theorem contains the classical theorems of Arzel\`{a}-Ascoli
and of Schauder. It can be considered as a special case, a ``lite'' version,
of considerably more abstract results presented by Robert G. Bartle in \cite
{bartle} (cf.\ also e.g., \cite{mujica}) and which have their roots in
earlier work of R.\ S.\ Phillips, \v {S}mulian and Kakutani. However it
seems easier to give a direct proof of this theorem than to deduce it from
\cite{bartle}.

\smallskip

\begin{theorem}
\label{aas}Let $A$ and $B$ be two sets and let $h:A\times B\rightarrow \Bbb{C%
}$ be a function with the properties that
\begin{equation}
\sup_{a\in A}\left| h(a,b)\right| <\infty \text{ for each fixed }b\in B\text{%
, and}  \label{qzi}
\end{equation}
\begin{equation}
\sup_{b\in B}\left| h(a,b)\right| <\infty \text{ for each fixed }a\in A.
\label{qzii}
\end{equation}

\noindent Define $d_{A}(a_{1},a_{2}):=\sup_{b\in B}\left|
h(a_{1},b)-h(a_{2},b)\right| $ for each pair of elements $a_{1}$ and $a_{2}$
in $A$.

\noindent Define $d_{B}(b_{1},b_{2})=\sup_{a\in A}\left|
h(a,b_{1})-h(a,b_{2})\right| $ for each pair of elements $b_{1}$ and $b_{2}$
in $B$.

Then
\begin{equation}
(A,d_{A})\text{ and }(B,d_{B})\text{ are semimetric spaces}  \label{qziii}
\end{equation}

\noindent and
\begin{equation}
(A,d_{A})\text{ is totally bounded if and only if }(B,d_{B})\text{ is
totally bounded. }  \label{qziv}
\end{equation}
\end{theorem}

\smallskip

\textit{Proof.} It is obvious that (\ref{qziii}) holds. For the proof of (%
\ref{qziv}), because of the symmetrical roles of $A$ and $B$, we only have
to prove one of the two implications.

Suppose then that $(A,d_{A})$ is totally bounded. By Proposition \ref{tbsep}%
, there exists a countable subset $Y$ of $A$ which is dense in $A$. Let us
show that
\begin{equation}
d_{B}(b_{1},b_{2})=\sup_{y\in Y}\left| h(y,b_{1})-h(y,b_{2})\right| \text{
for all }b_{1},b_{2}\in B.  \label{reyu}
\end{equation}
The inequality ``$\ge $'' in (\ref{reyu}) is obvious. For the reverse
inequality, given any $b_{1}$ and $b_{2}$ in $B$ and any arbitrarily small
positive $\epsilon $, we choose $a\in A$ such that
\begin{equation}
d_{B}(b_{1},b_{2})\le \left| h(a,b_{1})-h(a,b_{2})\right| +\epsilon /3.
\label{rzaq}
\end{equation}
Then we choose $z\in Y$ such that
\begin{equation}
d_{A}(z,a)<\epsilon /3.  \label{r3}
\end{equation}
It follows that \hskip 7.35 truemm $\displaystyle \left|
h(a,b_{1})-h(a,b_{2})\right|$

\vskip -5mm
\begin{eqnarray*}
\phantom{wwwwwssss}&\le &\left| h(a,b_{1})-h(z,b_{1})\right| +\left|
h(z,b_{1})-h(z,b_{2})\right| +\left| h(z,b_{2})-h(a,b_{2})\right| \\
&\le &2d_{A}(a,z)+\sup_{y\in Y}\left| h(y,b_{1})-h(y,b_{2})\right| .
\end{eqnarray*}
Combining this with (\ref{rzaq}) and (\ref{r3}), we can immediately complete
the proof of (\ref{reyu}).

\smallskip We shall now assume that $(B,d_{B})$ is not totally bounded and
show that this leads to a contradiction. By this assumption and by
Proposition \ref{phir}, there exists some positive number $r$ and some
infinite sequence $\left\{ b_{n}\right\} _{n\in \Bbb{N}}$ of elements of $B$
such that
\begin{equation}
d_{B}(b_{m},b_{n})>r\text{ for each }m,n\in \Bbb{N}\text{ with }m\ne n.
\label{rsep}
\end{equation}
For each fixed $y\in Y$ it follows from (\ref{qzii}) that the numerical
sequence $\left\{ h(y,b_{n})\right\} _{n\in \Bbb{N}}$ is bounded and thus
has a convergent subsequence. Since $Y$ is countable we can apply a standard
Cantor ``diagonalization'' argument to obtain a subsequence $\left\{
b_{\gamma _{n}}\right\} _{n\in \Bbb{N}}$ of $\left\{ b_{n}\right\} _{n\in
\Bbb{N}}$ such that $\lim_{n\rightarrow \infty }h(y,b_{\gamma _{n}})$ exists
for each $y\in Y$. Therefore, after simply changing our notation, we can
assume the existence of an infinite sequence $\left\{ b_{n}\right\} _{n\in
\Bbb{N}}$ in $B$ which satisfies (\ref{rsep}) and also
\begin{equation}
\lim_{n\rightarrow \infty }h(y,b_{n})\text{ exists and is finite for each }%
y\in Y.  \label{rlim}
\end{equation}

In view of (\ref{reyu}) and (\ref{rsep}), for each pair of integers $m$ and $%
n$ with $0<m<n$ there exists an element $y_{m,n}\in Y$ such that $\left|
h(y_{m,n},b_{m})-h(y_{m,n},b_{n})\right| >r$, and so, in particular,
\begin{equation}
\left| h(y_{m,m+1},b_{m})-h(y_{m,m+1},b_{m+1})\right| >r\text{ for all }m\in
\Bbb{N}.  \label{sqwe}
\end{equation}
Our assumption that $(A,d_{A})$ is totally bounded ensures, by Proposition
\ref{koshy}, that there exists a strictly increasing sequence of positive
integers $\left\{ m_{k}\right\} _{k\in \Bbb{N}}$ such that $\left\{
y_{m_{k},m_{k}+1}\right\} _{k\in \Bbb{N}}$ is a Cauchy sequence in $%
(A,d_{A}) $. Now we set $z_{k}=y_{m_{k},m_{k}+1}$ for each $k$. We choose
some sufficiently large integer $N$ for which
\begin{equation}
d_{A}(z_{N},z_{k})<r/4\text{ for all }k\ge N.  \label{wsx}
\end{equation}
Now we combine (\ref{sqwe}) and (\ref{wsx}) to obtain that, for each $k\ge N$%
,
\begin{eqnarray*}
r &<&\left| h(z_{k},b_{m_{k}})-h(z_{k},b_{m_{k}+1})\right| \\
&\le &\left| h(z_{k},b_{m_{k}})-h(z_{N},b_{m_{k}})\right| +\left|
h(z_{N},b_{m_{k}})-h(z_{N},b_{m_{k}+1})\right| \\
&&+\left| h(z_{N},b_{m_{k}+1})-h(z_{k},b_{m_{k}+1})\right| \\
&<&\frac{r}{4}+\left| h(z_{N},b_{m_{k}})-h(z_{N},b_{m_{k}+1})\right| +\frac{r%
}{4}.
\end{eqnarray*}
In view of (\ref{rlim}), we obtain that $\lim_{k\rightarrow \infty }\left|
h(z_{N},b_{m_{k}})-h(z_{N},b_{m_{k}+1})\right| =0$. So the inequalities on
the preceding lines would imply that $r\le r/2$.
This contradiction shows that $(B,d_{B})$ must be totally bounded, and so
completes the proof of the theorem. \qed

\smallskip

\section{\label{pll}Preliminaries about lattices and lattice couples.}

Let $(\Omega ,\Sigma ,\mu )$ be a $\sigma $-finite measure space here and in
the sequel. (Some of the assertions which we will be making here are simply
false if $(\Omega ,\Sigma ,\mu )$ is not $\sigma $-finite.)

\begin{definition}
\label{cbl}We say that a Banach space $X$ is a \textbf{\textit{CBL}}, or a
\textbf{\textit{complexified Banach lattice of measurable functions on}} $%
\Omega $ if

(i) all the elements of $X$ are (equivalence classes of a.e.\ equal)
measurable functions $f:\Omega \rightarrow \Bbb{C}$ and

(ii) for any measurable functions $f:\Omega \rightarrow \Bbb{C}$ and $%
g:\Omega \rightarrow \Bbb{C}$, if $f\in X$ and $\left| g\right| \le \left|
f\right| $ a.e., then $g\in X$ and $\left\| g\right\| _{X}\le \left\|
f\right\| _{X}$.
\end{definition}

We will now recall a number of definitions and basic facts about CBLs. In
several cases the relevant proofs of these facts in the literature to which
we refer are given for Banach lattices of \textit{real} valued functions.
But in all those cases it is an obvious and easy exercise to adapt those
proofs to our case here.

Any two CBLs $X_{0}$ and $X_{1}$ on the same underlying measure space always
form a Banach couple. See e.g., \cite{ca} p.\ 122 and p.\ 161 or \cite{kps}
Corollary 1, p.\ 42, or \cite{cn} Remark 1.41, pp.\ 34--35. (As explicitly
stated and shown in \cite{cn}, this is also true for non $\sigma$-finite
measure spaces.)

\smallskip For each CBL $X$ on $(\Omega ,\Sigma ,\mu )$, there exists a
measurable subset $\Omega _{X}$ of $\Omega $, which may be called the
\textbf{\textit{support}} of $X$, such that, for every function $g\in X$, we
have $g(\omega )=0$ for a.e.\ $\omega \in \Omega \backslash \Omega _{X}$.
Furthermore, there exists a function $f_{X}\in X$ such that $f_{X}(\omega
)>0 $ for a.e.\ $\omega \in \Omega _{X}$. (Cf.\ e.g., Remarks 1.3 and 1.4 on
p.\ 14 of \cite{cn}.) Obviously the set $\Omega _{X}$ is unique to within a
set of measure zero. (Of course, on the other hand, the function $f_{X}$
certainly is \textit{not} unique.) If $\Omega _{X}=\Omega $ (at least to
within a set of measure zero) then we say that $X$ is \textbf{\textit{%
saturated}}.

The set $\Omega _{X}$ has an additional useful property: There exists a
sequence of sets $\left\{ E_{n}\right\} _{n\in \Bbb{N}}$ in $\Sigma $ such
that
\begin{equation}
\Omega _{X}=\bigcup_{n\in \Bbb{N}}E_{n}\text{ with }E_{n}\subset E_{n+1},\mu
(E_{n})<\infty \text{, and }\chi _{E_{n}}\in X\text{ for each }n\in \Bbb{N}.
\label{yyii}
\end{equation}
The actual construction of $\Omega _{X}$ and of the sequence $\left\{
E_{n}\right\} _{n\in \Bbb{N}}$ can be performed by an ``exhaustion'' process
described in the proof of Theorem 3 on pp.\ 455--456 of \cite{z} and also
described (perhaps slightly more explicitly for our purposes here) in the
first part of the proof of Proposition 4.1 on p.\ 58 of \cite{cn1}. (Note
however that there is a small misprint in \cite{cn1}, the omission of ``$\mu
(E)$'', in the third line of this latter proof. I.e., the numbers $\alpha
_{k}$ must of course be defined by $\alpha _{k}=\sup \left\{ \mu (E):E\in
\Sigma ,E\subset F_{k},\chi _{E}\in X\right\} $.) For one possible (very
easy and of course not unique) way to construct a function $f_{X}\in X$ with
the above mentioned property see, e.g., \cite{cn} p.\ 14 Remark 1.4.

We will say that $(X_{0},X_{1})$ is a \textit{\textbf{saturated lattice
couple}}, if it is a lattice couple and both $X_{0}$ and $X_{1}$ are
saturated.

\begin{lemma}
\smallskip \label{bosat}If $(X_{0},X_{1})$ is a saturated lattice couple
then $X_{0}\cap X_{1}$ is saturated, and $[X_{0},X_{1}]_{\theta }$ is
saturated for each $\theta \in (0,1)$.
\end{lemma}

\textit{Proof. }Let $(\Omega ,\Sigma ,\mu )$ be the underlying measure space
for the couple. The function $\min \{f_{X_{0}},f_{X_{1}}\}$ is in $X_{0}\cap
X_{1}$ and therefore it is also in $[X_{0},X_{1}]_{\theta }$. It is strictly
positive a.e.\ on $\Omega $. So neither of the sets $\Omega \backslash
\Omega _{X_{0}\cap X_{1}}$ and $\Omega \backslash \Omega
_{[X_{0},X_{1}]_{\theta }}$ can have positive measure. \qed

\smallskip

Given an arbitrary CBL $X$ on $(\Omega ,\Sigma ,\mu )$ we define the
functional $\left\| \cdot \right\| _{X^{\prime }}$ by
\begin{equation}
\left\| f\right\| _{X^{\prime }}:=\sup \left\{ \left| \int_{\Omega }fgd\mu
\right| :g\in X,\left\| g\right\| _{X}\le 1\right\}  \label{czd}
\end{equation}
for each measurable function $f:\Omega \rightarrow \Bbb{C}$.

\begin{remark}
\label{ocr}Obviously we can replace $\left| \int_{\Omega }fgd\mu \right| $
by $\int_{\Omega }|fg|d\mu $ in the formula (\ref{czd}).
\end{remark}

Let $X^{\prime }$ be the set of all measurable functions $f:\Omega
\rightarrow \Bbb{C}$ for which $\left\| f\right\| _{X^{\prime }}<\infty $.
Clearly $X^{\prime }$ is a linear space and $\left\| \cdot \right\|
_{X^{\prime }}$ is a seminorm on $X^{\prime }$ satisfying
\begin{equation}
\left| \int_{\Omega }fgd\mu \right| \le \left\| f\right\| _{X^{\prime
}}\left\| g\right\| _{X}\text{ for all }f\in X^{\prime }\text{ and all }g\in
X.  \label{erzz}
\end{equation}
The space $X^{\prime }$ is customarily referred to as the \textbf{\textit{%
K\"{o}the dual}} or the \textbf{\textit{associate space}} of $X$.

If $\mu (\Omega _{X})>0$, then, via a series of theorems, including one (%
\cite{z} Theorem 1, p.\ 470) which uses Hilbert space techniques, it can be
shown that $X^{\prime }$ is non trivial, i.e., it contains elements which do
not vanish a.e.\ on $\Omega _{X}$. If, furthermore, $X$ is saturated, then $%
\left\| \cdot \right\| _{X^{\prime }}$ is a norm with respect to which $%
X^{\prime }$ is a saturated CBL on $(\Omega ,\Sigma ,\mu )$. (See e.g., \cite
{z} p.\ 472, Theorem 4.)

\smallskip Of course $X^{\prime }$ can be identified with a subspace of $%
X^{*}$, the dual space of $X$, and in some, but not all, cases it is also a
\textit{\textbf{norming subspace }}of $X^{*}$, i.e., it satisfies
\begin{equation}
\left\| g\right\| _{X}=\sup \left\{ \left| \int_{\Omega }fgd\mu \right|
:f\in X^{\prime },\left\| f\right\| _{X^{\prime }}\le 1\right\} \text{ for
each }g\in X.  \label{ains}
\end{equation}

\smallskip The following result of Lorentz and Luxemburg, which appears as
Proposition 1.b.18 on p.\ 29 of \cite{lt}, gives necessary and sufficient
conditions on $X$ for (\ref{ains}) to hold. In particular, it implies that
the $\sigma $-order continuity of $X$ is a sufficient condition. So is the
Fatou property.

\begin{theorem}
\label{lolu}Let $X$ be an arbitrary CBL on a $\sigma $-finite measure space $%
(\Omega ,\Sigma ,\mu )$. The associate space $X^{\prime }$ is a norming
subspace of $X^{*}$ if and only if $\lim_{n\rightarrow \infty }\left\|
f_{n}\right\| _{X}=\left\| f\right\| _{X}$ for every non negative function $%
f\in X$ and every sequence $\left\{ f_{n}\right\} _{n\in \Bbb{N}}$ of
measurable functions satisfying $0\le f_{n}(\omega )\le f_{n+1}(\omega )\le
f(\omega )$ for a.e.\ $\omega $ and each $n\in \Bbb{N}$ and $%
\lim_{n\rightarrow \infty }f_{n}(\omega )=f(\omega )$ for a.e.\ $\omega $.
\end{theorem}

\textit{Proof. }The proof is essentially the same as the proof of
Proposition 1.b.18 on p.\ 29 of \cite{lt}. Some small modifications are
required in lines 6 and 7 of that proof on p.\ 30 of \cite{lt}, because in
the statement and proof of that proposition it is assumed that $X$ is a
K\"{o}the function space. Here is the required replacement for those two
lines, using the notation of \cite{lt}:

Let now $f\in X$ be an element with $\left\| f\right\| >1$. Using Remark 1.3
on p.\ 14 of \cite{cn}, we know that there exists an expanding sequence $%
\left\{ E_{k}\right\} _{k\in \Bbb{N}}$ of measurable subsets of $\Omega $
such that $\chi _{E_{k}}\in X$ and $\mu (E_{k})<\infty $ for each $k$, and
such that $f$ vanishes a.e.\ on the complement of $\bigcup_{k\in \Bbb{N}%
}E_{k}$. Then $\left| f\right| \chi _{E_{k}}\uparrow \left| f\right| $ a.e.\
and so $\left\| f\chi _{\sigma }\right\| >1$ for some set $\sigma =E_{k}$
for some sufficiently large $k$. By using the separation theorem for $Y$.....

From here the proof can be continued and concluded exactly as in \cite{lt}.
\qed

The associate space $(X^{\prime })^{\prime }$ of $X^{\prime }$, i.e.\ the
\textit{second associate} of $X$ is usually denoted by $X^{\prime \prime }$.
Obviously $X\subset X^{\prime \prime }$ and $\left\| x\right\| _{X^{\prime
\prime }}\le \left\| x\right\| _{X}$ for each $x\in X$. Obviously $X^{\prime
\prime }$ is a CBL whenever $X$ (and therefore also $X^{\prime }$) is
saturated.

\smallskip As in e.g., \cite{kps}, we say that the CBL $X$ has \textbf{%
\textit{absolutely continuous norm}} if $\lim_{n\rightarrow \infty }\left\|
f\chi _{E_{n}}\right\| _{X}=0$ for every $f\in X$ and every sequence $%
\left\{ E_{n}\right\} _{n\in \Bbb{N}}$ of measurable sets satisfying $%
E_{n+1}\subset E_{n}$ for all $n$ and $\bigcap_{n\in \Bbb{N}}E_{n}=\emptyset
$. As in e.g., \cite{lt}, we say that the CBL $X$ is $\sigma $\textbf{-%
\textit{order continuous} }if $\lim_{n\rightarrow \infty }\left\|
f_{n}\right\| _{X}=0$ for every sequence $\left\{ f_{n}\right\} _{n\in \Bbb{N%
}}$ of functions in $X$ satisfying $0\le f_{n+1}\le f_{n}$ and $%
\lim_{n\rightarrow \infty }f_{n}=0$ a.e. It is easy to see that these two
properties of $X$ are in fact equivalent.

A CBL $X$ is said to have the \textbf{\textit{Fatou property}} if whenever $%
\left\{ f_{n}\right\} _{n\in \Bbb{N}}$ is a norm bounded a.e.\ monotonically
non decreasing sequence of non negative functions in $X$, its a.e.\
pointwise limit $f$ is also in $X$ with $\left\| f\right\|
_{X}=\lim_{n\rightarrow \infty }\left\| f_{n}\right\| _{X}$. If $X$ is
saturated, then $X$ has the Fatou property if and only if $X=X^{\prime
\prime }$ isometrically. (See \cite{z} p.\ 472. Cf.\ also \cite{lt} p.\ 30,
but recall that there extra hypotheses are imposed.)

\smallskip We remark that obvious counterexamples (see e.g., \cite{cn}
Remark 7.3 p.\ 92) show that the above claims about $X^{\prime }$ and $%
X^{\prime \prime }$ are false for certain non $\sigma $-finite measure
spaces.

Given a pair of CBLs $X_{0}$ and $X_{1}$ on $(\Omega ,\Sigma ,\mu )$ and a
number $\theta \in (0,1)$, we define the space $X_{0}^{1-\theta
}X_{1}^{\theta }$ , analogously to the definition in \cite{ca} Section 13.5
pp.\ 123, to be the set of all measurable functions $f:\Omega \rightarrow
\Bbb{C}$ of the form
\begin{equation}
f=uf_{0}^{1-\theta }f_{1}^{\theta }  \label{splrc}
\end{equation}
where $u\in L^{\infty }(\mu )$ and $f_{j}$ is a non negative function in $%
\mathcal{B}_{X_{j}}$ for $j=0,1$. For each $f\in X_{0}^{1-\theta
}X_{1}^{\theta }$ we define $\left\| f\right\| _{X_{0}^{1-\theta
}X_{1}^{\theta }}=\inf \left\| u\right\| _{L^{\infty }(\mu )}$, where the
infimum is taken over all representations of $f$ of the form (\ref{splrc})
with the stated properties. It can be shown that this is in fact a norm on $%
X_{0}^{1-\theta }X_{1}^{\theta }$, with respect to which $X_{0}^{1-\theta
}X_{1}^{\theta }$ is a CBL. This is proved in Section 33.5 on pp.\ 164--165
of \cite{ca}.

The norm $1$ inclusions
\begin{equation}
\lbrack X_{0},X_{1}]_{\theta }\overset{1}{\subset }X_{0}^{1-\theta
}X_{1}^{\theta }\overset{1}{\subset }[X_{0},X_{1}]^{\theta }  \label{xsw}
\end{equation}
are special cases (set $B_{0}=B_{1}=\Bbb{C}$) of the results (i) and (ii) of
Section 13.6 on p.\ 125\ of \cite{ca} (proved in \cite{ca} Section 33.6 on
pp.\ 171--180). Furthermore, with the help of Bergh's theorem \cite{bergh}, (%
\ref{xsw}) can be strengthened to tell us that
\begin{equation}
\left\| x\right\| _{[X_{0},X_{1}]_{\theta }}=\left\| x\right\|
_{X_{0}^{1-\theta }X_{1}^{\theta }}=\left\| x\right\|
_{[X_{0},X_{1}]^{\theta }}\text{ for all }x\in [X_{0},X_{1}]_{\theta }
\label{btu}
\end{equation}

We will need to use the formula
\begin{equation}
\left( X_{0}^{1-\theta }X_{1}^{\theta }\right) ^{\prime }=(X_{0}^{\prime
})^{1-\theta }(X_{1}^{\prime })^{\theta }  \label{loz}
\end{equation}
which holds with equality of norms (or seminorms when $\Omega _{X_{0}}$ or $%
\Omega _{X_{1}}$ is strictly smaller than $\Omega $) for all pairs of CBLs $%
X_{0}$ and $X_{1}$ on $(\Omega ,\Sigma ,\mu )$. This formula was originally
stated and proved by Lozanovskii under certain hypotheses, then by Reisner
\cite{reisner} under other hypotheses. The general version stated here is
proved in \cite{cn} Section 7, pp.\ 91--97 using Reisner's proof and a
remark of Kalton.

We will also need this lemma:

\begin{lemma}
\smallskip \label{fca}Let $(X_{0},X_{1})$ be an arbitrary saturated lattice
couple. For each $\theta \in (0,1)$ we have
\end{lemma}

\begin{equation}
\sup_{y\in \mathcal{B}_{[X_{0}^{\prime },X_{1}^{\prime }]_{\theta }}}\left|
\int_{\Omega }xyd\mu \right| =\sup_{y\in \mathcal{B}_{(X_{0}^{\prime
})^{1-\theta }(X_{1}^{\prime })^{\theta }}}\left| \int_{\Omega }xyd\mu
\right| \text{ for each }x\in X_{0}\cap X_{1}.  \label{ffza}
\end{equation}

\textit{Proof. }Applying (\ref{btu}) to the couple $(X_{0}^{\prime
},X_{1}^{\prime })$, we of course obtain the inequality $``\le "$ in (\ref
{ffza}). To show the reverse inequality $``\ge "$, we fix some $x\in
X_{0}\cap X_{1}$ and $y\in \mathcal{B}_{(X_{0}^{\prime })^{1-\theta
}(X_{1}^{\prime })^{\theta }}$ and we shall construct a sequence $\left\{
y_{n}\right\} _{n\in \Bbb{N}}$ in $\mathcal{B}_{[X_{0}^{\prime
},X_{1}^{\prime }]_{\theta }}$ for which
\begin{equation}
\lim_{n\rightarrow \infty }\int_{\Omega }xy_{n}d\mu =\int_{\Omega }xyd\mu .
\label{wcya}
\end{equation}
By Lemma \ref{bosat}, since $X_{0}^{\prime }$ and $X_{1}^{\prime }$ are both
saturated, so is $[X_{0}^{\prime },X_{1}^{\prime }]_{\theta }$. Consequently
(cf.\ (\ref{yyii})) there exists an expanding sequence $\left\{
E_{n}\right\} _{n\in \Bbb{N}}$ of sets in $\Sigma $ such that $\Omega
=\bigcup_{n\in \Bbb{N}}E_{n}$ and $\chi _{E_{n}}\in [X_{0}^{\prime
},X_{1}^{\prime }]_{\theta }$ for each $n\in \Bbb{N}$. Let $y_{n}=y\chi
_{E_{n}\cap \{\omega \in \Omega :\left| y(\omega )\right| \le n\}}$. Then $%
y_{n}\in [X_{0}^{\prime },X_{1}^{\prime }]_{\theta }$ and we have $\left|
xy_{n}\right| \le \left| xy\right| $ and $\lim_{n\rightarrow \infty
}x(\omega )y_{n}(\omega )=x(\omega )y(\omega )$ for all $\omega \in \Omega $%
. The function $xy$ is integrable, since $X_{0}\cap X_{1}\overset{1}{\subset
}[X_{0},X_{1}]_{\theta }\overset{1}{\subset }X_{0}^{1-\theta }X_{1}^{\theta
} $ and $(X_{0}^{\prime })^{1-\theta }(X_{1}^{\prime })^{\theta }\overset{1}{%
=}\left( X_{0}^{1-\theta }X_{1}^{\theta }\right) ^{\prime }$ (cf.\ (\ref{btu}%
) and (\ref{loz}) and Remark \ref{ocr}). So (\ref{wcya}) follows from the
Lebesgue dominated convergence theorem. This completes the proof. \qed

\section{\label{main}The main result.}

\smallskip Our main result is a corollary of the following theorem.

\begin{theorem}
\smallskip \label{rlat} Let $\vec{G}=(G_{0},G_{1})$ be an arbitrary Banach
couple and let $\vec{X}=(X_{0},X_{1})$ be a saturated lattice couple. Then
every linear operator $T$ which satisfies $T:\vec{G}\overset{c,b}{%
\rightarrow }\vec{X}$ has the compactness property
\begin{equation}
T:[G_{0},G_{1}]_{\theta }\overset{c}{\rightarrow }[X_{0}^{\prime \prime
},X_{1}^{\prime \prime }]_{\theta }  \label{hcp}
\end{equation}
for each $\theta \in (0,1)$.
\end{theorem}

\begin{corollary}
\label{crlat}Let $\vec{X}=(X_{0},X_{1})$ be a saturated lattice couple.
Suppose that either

(i) $X_{0}$ and $X_{1}$ both have the Fatou property, or

(ii) at least one of the spaces $X_{0}$ and $X_{1}$ is $\sigma $-order
continuous, or

(iii) there exists at least one value of $\theta \in (0,1)$ for which $%
X=X_{0}^{1-\theta }X_{1}^{\theta }$ satisfies the condition of Theorem \ref
{lolu}.

Then
\begin{equation*}
(*.*)\blacktriangleright \vec{X}.
\end{equation*}
\end{corollary}

\begin{remark}
\label{mnotsat} The requirement that $(X_{0},X_{1})$ is saturated is merely
a technical convenience which makes the formulation and proof of Theorem \ref
{rlat} simpler and shorter. In fact it is entirely unnecessary for Corollary
\ref{crlat}. The easy and rather obvious argument which extends the proof of
Corollary \ref{crlat} to the non saturated case uses the easily checked fact
that $\Omega _{X_{0}\cap X_{1}}=\Omega _{[X_{0},X_{1}]_{\theta }}$ and
replaces the spaces $X_{0}$ and $X_{1}$ in an appropriate way by their
``restrictions'' to the smaller measure space $\Omega _{X_{0}\cap X_{1}}$.
The details of that argument are given below in Subsection \ref{dnotsat}.
\end{remark}

\textit{Proof of Theorem \ref{rlat}.} Since $[G_{0}^{\circ },G_{1}^{\circ
}]_{\theta }=[G_{0},G_{1}]_{\theta }$ (\cite{ca} Sections 9.3 (p.\ 116) and
29.3 (pp.\ 113--4)) we can clearly suppose without loss of generality that $%
\vec{G}$ is a regular couple. Let $\left\langle \cdot ,\cdot \right\rangle $
denote the duality between $G_{0}\cap G_{1}$ and $\left( G_{0}\cap
G_{1}\right) ^{*}$. Let $G$ be any one of the spaces $G_{0}$, $G_{1}$ or $%
[G_{0},G_{1}]_{\theta }$ and define $G^{\#}$ to be the subspace of elements $%
\gamma \in \left( G_{0}\cap G_{1}\right) ^{*}$ for which the norm $\left\|
\gamma \right\| _{G^{\#}}:=$ $\sup \left\{ \left| \left\langle g,\gamma
\right\rangle \right| :g\in \mathcal{B}_{G}\cap G_{0}\cap G_{1}\right\} $ is
finite. Of course $G^{\#}$, when equipped with this norm, is a Banach space
which is continuously embedded in $\left( G_{0}\cap G_{1}\right) ^{*}$. So $%
(G_{0}^{\#},G_{1}^{\#})$ is a Banach couple.

We could of course identify $G^{\#}$ with the dual of $G$, but it is more
convenient to use the above definition. Note also that in fact $%
G_{0}^{\#}+G_{1}^{\#}\overset{1}{=}\left( G_{0}\cap G_{1}\right) ^{*}$.
Calder\'{o}n's remarkable duality theorem (\cite{ca} Section 12.1 p.\ 121
and Section 32.1 pp.\ 148--156) can be expressed by the formula $\left(
[G_{0},G_{1}]_{\theta }\right) ^{\#}\overset{1}{=}[G_{0}^{\#},G_{1}^{\#}]^{%
\theta }$. For a more detailed discussion of all these issues we refer to
\cite{dualitynotes}.

Let $T$ be an arbitrary linear operator satisfying $T:\vec{G}\overset{c,b}{%
\rightarrow }\vec{X}$. We may suppose, without loss of generality, that $%
\left\| T\right\| _{\vec{G}\rightarrow \vec{X}}:=\max_{j=0,1}\left\|
T\right\| _{G_{j}\rightarrow X_{j}}=1$. For $j=0,1$, let $X_{j}^{\prime }$
be the associate space of $X_{j}$. Let $(\Omega ,\Sigma ,\mu )$ be the
underlying measure space for $(X_{0},X_{1})$. For each $g\in G_{0}\cap G_{1}$
and each $z\in X_{0}^{\prime }+X_{1}^{\prime }$ define $h(g,z)=\int_{\Omega
}zTgd\mu $. Of course (cf.\ (\ref{erzz})) the function $h$ satisfies
\begin{equation}
\left| h(g,z)\right| \le \left\| z\right\| _{X_{j}^{\prime }}\left\|
Tg\right\| _{X_{j}}\le \left\| z\right\| _{X_{j}^{\prime }}\left\| g\right\|
_{G_{j}}  \label{ercv}
\end{equation}
for $j=0,1$ and all $g\in G_{0}\cap G_{1}$ and $z\in X_{j}^{\prime }$.
Therefore $h$ also satisfies
\begin{equation}
\left| h(g,z)\right| \le \left\| z\right\| _{X_{0}^{\prime }+X_{1}^{\prime
}}\left\| g\right\| _{G_{0}\cap G_{1}}\text{ for all }g\in G_{0}\cap G_{1}%
\text{ and }z\in X_{0}^{\prime }+X_{1}^{\prime }.  \label{tcfm}
\end{equation}

\smallskip For each fixed $z\in X_{0}^{\prime }+X_{1}^{\prime }$ we define
the linear functional $Sz$ on $G_{0}\cap G_{1}$ by $\left\langle
g,Sz\right\rangle =h(g,z)$. Of course $Sz$ depends linearly on $z$ and it is
clear from (\ref{tcfm}) that we have thus defined a bounded linear operator $%
S:X_{0}^{\prime }+X_{1}^{\prime }\rightarrow \left( G_{0}\cap G_{1}\right)
^{*}$. For $j=0,1$, in view of (\ref{ercv}), we see that, for each $z\in
X_{j}^{\prime }$, we have $Sz\in G_{j}^{\#}$ with $\left\| Sz\right\|
_{G_{j}^{\#}}\le \left\| z\right\| _{X_{j}^{\prime }}$, i.e., $%
S:X_{j}^{\prime }\overset{b}{\rightarrow }G_{j}^{\#}.$ (Note, cf.\ \cite
{dualitynotes}, that we do not have to consider the extension of $Sz$ to a
space larger than $G_{0}\cap G_{1}$.)

We now wish to show that $S$ satisfies the compactness condition
\begin{equation}
S:X_{0}^{\prime }\overset{c}{\rightarrow }G_{0}^{\#}.  \label{hhyo}
\end{equation}
We will do this by applying Theorem \ref{aas}. We consider the restriction
of the function $h(g,y)$ to the set $A\times B$ where $A=\mathcal{B}%
_{G_{0}}\cap G_{1}$ and $B=\mathcal{B}_{X_{0}^{\prime }}$. Given any
sequence $\left\{ g_{n}\right\} _{n\in \Bbb{N}}$ in $A$, we of course have
(cf.\ (\ref{erzz}) and (\ref{ercv})) that $d_{A}(g_{m},g_{n})=\sup \left\{
\left| h(g_{m},z)-h(g_{n},z)\right| :z\in B\right\} \le \left\|
Tg_{m}-Tg_{n}\right\| _{X_{0}}$. So the fact that $T:G_{0}\overset{c}{%
\rightarrow }X_{0}$ implies that $(A,d_{A})$ is totally bounded. (Cf., e.g.,
Proposition \ref{phir} or Theorem 15 on p.\ 22 of \cite{dunfordschwartz}.)
Consequently, in view of Theorem \ref{aas} and Proposition \ref{koshy}, if $%
\left\{ z_{n}\right\} _{n\in \Bbb{N}}$ is an arbitrary sequence in $B$, then
it has a subsequence which is Cauchy with respect to the semimetric
\begin{eqnarray*}
d_{B}(y,z) &=&\sup \left\{ \left| h(g,y)-h(g,z)\right| :g\in A\right\} \\
&=&\sup \left\{ \left| \left\langle g,S(y-z)\right\rangle \right| :g\in
G_{0}\cap G_{1},\left\| g\right\| _{G_{0}}\le 1\right\} \\
&=&\left\| S(y-z)\right\| _{G_{0}^{\#}}
\end{eqnarray*}
This is exactly the condition (\ref{hhyo}).

Since $X_{0}^{\prime }$ and $X_{1}^{\prime }$ are both CBLs of measurable
functions on the measure space $(\Omega ,\Sigma ,\mu )$, we can use (\ref
{hhyo}) and $S:X_{1}^{\prime }\overset{b}{\rightarrow }G_{1}^{\#}$ and apply
part (c) of Corollary 7 on p.\ 270 of \cite{CwK} to deduce that
\begin{equation}
S:[X_{0}^{\prime },X_{1}^{\prime }]_{\theta }\overset{c}{\rightarrow }%
[G_{0}^{\#},G_{1}^{\#}]_{\theta }.  \label{hhqq}
\end{equation}

We are now ready for a second application of Theorem \ref{aas}. Once more we
will use the same function $h$ defined above and restricted to a set $%
A\times B$, where this time we choose $A=\mathcal{B}_{[G_{0},G_{1}]_{\theta
}}\cap G_{0}\cap G_{1}$ and $B=\mathcal{B}_{[X_{0}^{\prime },X_{1}^{\prime
}]_{\theta }}$. This time, for each $y,z\in B$, we of course have $S(y-z)\in
[G_{0}^{\#},G_{1}^{\#}]_{\theta }$. So, using the isometry $\left(
[G_{0},G_{1}]_{\theta }\right) ^{\#}\overset{1}{=}[G_{0}^{\#},G_{1}^{\#}]^{%
\theta }$ mentioned above, and then Bergh's theorem \cite{bergh}, we obtain
that
\begin{eqnarray*}
d_{B}(y,z) &=&\sup \left\{ \left| \left\langle g,S(y-z)\right\rangle \right|
:g\in \mathcal{B}_{[G_{0},G_{1}]_{\theta }}\cap G_{0}\cap G_{1}\right\}
=\left\| S(y-z)\right\| _{\left( [G_{0},G_{1}]_{\theta }\right) ^{\#}} \\
&=&\left\| S(y-z)\right\| _{[G_{0}^{\#},G_{1}^{\#}]^{\theta }}=\left\|
S(y-z)\right\| _{[G_{0}^{\#},G_{1}^{\#}]_{\theta }}.
\end{eqnarray*}
The compactness property (\ref{hhqq}) of $S$ implies that $(B,d_{B})$ is
totally bounded. Consequently, by Theorem \ref{aas}, $(A,d_{A})$ is also
totally bounded. In view of Proposition \ref{koshy} and the fact that $%
G_{0}\cap G_{1}$ is dense in $[G_{0},G_{1}]_{\theta }$ (\cite{ca} Section
9.3 (p.\ 116) and Section 29.3 (pp.\ 113--4)), this means that the proof of
Theorem \ref{rlat} will be complete once we have shown that
\begin{equation}
d_{A}(g_{1},g_{2})=\left\| Tg_{1}-Tg_{2}\right\| _{[X_{0}^{\prime \prime
},X_{1}^{\prime \prime }]_{\theta }}\text{ for all }g_{1},g_{2}\in A.
\label{ptw}
\end{equation}

By definition, for each $g_{1}$ and $g_{2}$ in $A$ we have
\begin{equation*}
d_{A}(g_{1},g_{2})=\sup_{y\in \mathcal{B}_{[X_{0}^{\prime },X_{1}^{\prime
}]_{\theta }}}\left| \int_{\Omega }yT(g_{1}-g_{2})d\mu \right|
\end{equation*}

At this stage we do not need to consider the particular form of the element $%
Tg_{1}-Tg_{2}$. We know that it is an element of $X_{0}\cap X_{1}$. So, to
obtain (\ref{ptw}) it suffices to show that
\begin{equation}
\sup_{y\in \mathcal{B}_{[X_{0}^{\prime },X_{1}^{\prime }]_{\theta }}}\left|
\int_{\Omega }xyd\mu \right| =\left\| x\right\| _{[X_{0}^{\prime \prime
},X_{1}^{\prime \prime }]_{\theta }}\text{ for each }x\in X_{0}\cap X_{1}.
\label{gllp}
\end{equation}

Since $X_{0}\cap X_{1}\subset X_{0}^{\prime \prime }\cap X_{1}^{\prime
\prime }\subset [X_{0}^{\prime \prime },X_{1}^{\prime \prime }]_{\theta }$,
we have, from (\ref{btu}) applied to the couple $(X_{0}^{\prime \prime
},X_{1}^{\prime \prime })$, that the right side of (\ref{gllp}) equals $%
\left\| x\right\| _{(X_{0}^{\prime \prime })^{1-\theta }(X_{1}^{\prime
\prime })^{\theta }}$ and this in turn, in view of (\ref{loz}) applied to
the couple $(X_{0}^{\prime },X_{1}^{\prime })$, equals $\left\| x\right\|
_{\left( (X_{0}^{\prime })^{1-\theta }(X_{1}^{\prime })^{\theta }\right)
^{\prime }}=\sup_{y\in \mathcal{B}_{(X_{0}^{\prime })^{1-\theta
}(X_{1}^{\prime })^{\theta }}}\left| \int_{\Omega }xyd\mu \right| $. Thus we
can complete the proof by applying Lemma \ref{fca}. \qed

\textit{Proof of Corollary \ref{crlat}. }As in the preceding proof, we
consider an arbitrary regular couple $(G_{0},G_{1})$ and an arbitrary
operator $T:\vec{G}\overset{c,b}{\rightarrow }\vec{X}$. We need to show that
\begin{equation}
T:[G_{0},G_{1}]_{\theta }\overset{c}{\rightarrow }[X_{0},X_{1}]_{\theta }%
\text{ }  \label{rtaa}
\end{equation}
for all $\theta \in (0,1)$. If condition (i) holds then $X_{0}^{\prime
\prime }=X_{1}$ and $X_{1}^{\prime \prime }=X_{1}$ and (\ref{hcp}) gives us
the required conclusion. Otherwise we simply work through all the same steps
as in the preceding proof until we reach (\ref{ptw}). Then (\ref{rtaa}) will
follow if, instead of (\ref{ptw}), we can establish a variant of (\ref{ptw})
or of (\ref{gllp}), namely that $\sup_{y\in \mathcal{B}_{[X_{0}^{\prime
},X_{1}^{\prime }]_{\theta }}}\left| \int_{\Omega }xyd\mu \right| =\left\|
x\right\| _{[X_{0},X_{1}]_{\theta }}$ for each $x\in X_{0}\cap X_{1}$. Using
Lemma \ref{fca} and then (\ref{loz}) and also (\ref{btu}), we see that the
required condition is equivalent to
\begin{equation}
\sup_{y\in \mathcal{B}_{\left( X_{0}^{1-\theta }X_{1}^{\theta }\right)
^{\prime }}}\left| \int_{\Omega }xyd\mu \right| =\left\| x\right\|
_{X_{0}^{1-\theta }X_{1}^{\theta }}\text{ for each }x\in X_{0}\cap X_{1}.
\label{ghkk}
\end{equation}
In fact condition (ii) implies condition (iii) because, if $X_{0}$ or $X_{1}$
is $\sigma $-order continuous, then $X_{0}^{1-\theta }X_{1}^{\theta }$ is $%
\sigma $-order continuous for every $\theta \in (0,1)$ (cf.\ Proposition 4
on p.\ 80 of \cite{reisner} or Theorem 1.29 on p.\ 27 of \cite{cn}).
Condition (iii) ensures that $\left( X_{0}^{1-\theta }X_{1}^{\theta }\right)
^{\prime }$ is a norming subspace of $\left( X_{0}^{1-\theta }X_{1}^{\theta
}\right) ^{*}$ for at least one value of $\theta $. This implies that (\ref
{ghkk}) holds for that value of $\theta $. Consequently, (\ref{rtaa}) holds
for that same value of $\theta $. So, by ``extrapolation'' (see Theorem 2.1,
on p.\ 339 of \cite{duke} or Theorem 5.3 on p.\ 311 of \cite{cks}), we
obtain (\ref{rtaa}) for \textit{all} $\theta \in (0,1)$. \qed

\smallskip

\section{\label{fpg}Further possible generalizations.}

\smallskip

\subsection{\label{dnotsat}\protect\smallskip Extending the result to the
case where $(X_{0},X_{1})$ is not saturated.}

Here we give a more detailed (perhaps too detailed?) explanation of the
claim made in Remark \ref{mnotsat}.

\begin{corollary}
\smallskip \label{zzlat}The result of Corollary \ref{crlat} also holds if $%
X_{0}$ is not saturated and/or $X_{1}$ is not saturated.
\end{corollary}

\textit{Proof. }Since $\Omega _{X_{0}\cap X_{1}}=\left\{ \omega \in \Omega
:f_{X_{0}\cap X_{1}}(\omega )>0\right\} $ for some non negative function $%
f_{X_{0}\cap X_{1}}$ which is in $X_{0}\cap X_{1}$ and therefore also in $%
[X_{0},X_{1}]_{\theta }$, we see that
\begin{equation*}
\mu \left( \Omega _{X_{0}\cap X_{1}}\backslash \Omega
_{[X_{0},X_{1}]_{\theta }}\right) =\mu \left( \left\{ \omega \in \Omega
:f_{X_{0}\cap X_{1}}(\omega )>0,\omega \notin \Omega _{[X_{0},X_{1}]_{\theta
}}\right\} \right) =0.
\end{equation*}
Next we remark that $\Omega _{[X_{0},X_{1}]_{\theta }}=\left\{ \omega \in
\Omega :f_{[X_{0},X_{1}]_{\theta }}(\omega )>0\right\} $ for some non
negative function $f_{[X_{0},X_{1}]_{\theta }}$ which is in $%
[X_{0},X_{1}]_{\theta }$ and therefore in $X_{0}^{1-\theta }X_{1}^{\theta }$%
. We thus have $f_{[X_{0},X_{1}]_{\theta }}=f_{0}^{1-\theta }f_{1}^{\theta }$
where $f_{0}$ and $f_{1}$ are non negative functions in $X_{0}$ and $X_{1}$
respectively. Clearly
\begin{equation}
v:=\min \left\{ f_{0},f_{1}\right\}  \label{inv}
\end{equation}
is a non negative function in $X_{0}\cap X_{1}$ and
\begin{eqnarray*}
\mu \left( \Omega _{[X_{0},X_{1}]_{\theta }}\backslash \Omega _{X_{0}\cap
X_{1}}\right) &=&\mu \left( \left\{ \omega \in \Omega
:f_{[X_{0},X_{1}]_{\theta }}(\omega )>0,\omega \notin \Omega _{X_{0}\cap
X_{1}}\right\} \right) \\
&=&\mu \left( \left\{ \omega \in \Omega :v(\omega )>0,\omega \notin \Omega
_{X_{0}\cap X_{1}}\right\} \right) =0.
\end{eqnarray*}
Thus we have shown that
\begin{equation}
\Omega _{X_{0}\cap X_{1}}=\Omega _{[X_{0},X_{1}]_{\theta }}\text{ a.e. }
\label{pnt}
\end{equation}

Now we consider the measure space $\left( \Omega _{*},\Sigma _{*},\mu
_{*}\right) $ where $\Omega _{*}=\Omega _{X_{0}\cap X_{1}}$ and $\Sigma _{*}$
is the $\sigma $-algebra of all sets in $\Sigma $ which are contained in $%
\Omega _{*}$, and $\mu _{*}$ is the restriction of $\mu $ to $\Sigma _{*}$.
We shall conveniently ``navigate'' between spaces of functions on $\Omega $
and spaces of function on $\Omega _{*}$ with the help of two simple and
obvious operators $\mathcal{R}$ and $\mathcal{E}$ of restriction and
extension. For each function $f:\Omega \rightarrow \Bbb{C}$ let $\mathcal{R}%
f $ be the restriction of $f$ to $\Omega _{*}$. For each function $g:\Omega
_{*}\rightarrow \Bbb{C}$ let $\mathcal{E}g$ be the complex valued function
on $\Omega $ which equals $0$ on $\Omega \backslash \Omega _{*}$ and
coincides with $g$ on $\Omega _{*}$. For $j=0,1$ we let $Y_{j}=\mathcal{R}%
X_{j}$. Thus $Y_{j}$ is a space of $\mu _{*}$ measurable functions $y:\Omega
_{*}\rightarrow \Bbb{C}$ and we may norm it by setting $\left\| y\right\|
_{Y_{j}}=\left\| \mathcal{E}y\right\| _{X_{j}}$. It is clear that $Y_{j}$ is
a CBL. Furthermore it is saturated, because the function $\mathcal{R}v$
(where $v$ is the function introduced in (\ref{inv})) is in $Y_{j}$ and is
strictly positive a.e.\ on $\Omega _{*}$. Obviously $\mathcal{E}:Y_{j}%
\overset{b}{\rightarrow }X_{j}$ and $\mathcal{R}:X_{j}\overset{b}{%
\rightarrow }Y_{j}$ for $j=0,1$ with norm $1$ in each case. It follows by
interpolation that $\mathcal{E}:[Y_{0},Y_{1}]_{\theta }\overset{b}{%
\rightarrow }[X_{0},X_{1}]_{\theta }$ and $\mathcal{R}:[X_{0},X_{1}]_{\theta
}\overset{b}{\rightarrow }[Y_{0},Y_{1}]_{\theta }$ (also in fact with norm $%
1 $).

Suppose now that given couple $(X_{0},X_{1})$ satisfies one (or more) of the
conditions (i), (ii) and (iii) stated in Corollary \ref{crlat}. Then the
couple $(Y_{0},Y_{1})$ satisfies the same condition: If $X_{j}$ has the
Fatou property then so, obviously does $Y_{j}$. Also, if $X_{j}$ is $\sigma $%
-order continuous, so is $Y_{j}$. To deal with condition (iii) we first
remark that it is easy to see that $\left\| g\right\| _{Y_{0}^{1-\theta
}Y_{1}^{\theta }}=\left\| \mathcal{E}g\right\| _{X_{0}^{1-\theta
}X_{1}^{\theta }}$ for each $g\in Y_{0}^{1-\theta }Y_{1}^{\theta }$.
Consequently, if $X_{0}^{1-\theta }X_{1}^{\theta }$ satisfies the hypotheses
of Theorem \ref{lolu} for some value of $\theta $, then so does $%
Y_{0}^{1-\theta }Y_{1}^{\theta }$.\smallskip

\smallskip Let $\vec{G}$ be an arbitrary Banach couple and suppose that $T:%
\vec{G}\overset{c,b}{\rightarrow }\vec{X}$. Then the composed operator $%
\mathcal{R}T$ satisfies $\mathcal{R}T:\vec{G}\overset{c,b}{\rightarrow }\vec{%
Y}$. Thus, by Corollary \ref{crlat}, we have $\mathcal{R}T:[G_{0},G_{1}]_{%
\theta }\overset{c}{\rightarrow }[Y_{0},Y_{1}]_{\theta }$. Consequently
\begin{equation}
\mathcal{ER}T:[G_{0},G_{1}]_{\theta }\overset{c}{\rightarrow }%
[X_{0},X_{1}]_{\theta }.  \label{bonk}
\end{equation}
In view of (\ref{pnt}), each function $f\in [X_{0},X_{1}]_{\theta }$
vanishes a.e.\ on $\Omega \backslash \Omega _{*}$, and therefore satisfies $%
\mathcal{ER}f=f$. This means that $\mathcal{ER}Tg=Tg$ for each $g\in
[G_{0},G_{1}]_{\theta }$. So (\ref{bonk}) gives us that $T:[G_{0},G_{1}]_{%
\theta }\overset{c}{\rightarrow }[X_{0},X_{1}]_{\theta }$ and completes the
proof of Corollary \ref{zzlat}. \qed

\subsection{Two sided interpolation of compactness.}

We may use the notation $T:\vec{A}\overset{c,c}{\rightarrow }\vec{B}$ to
mean that the linear operator $T:A_{0}+A_{1}\rightarrow B_{0}+B_{1}$
satisfies the ``two sided'' compactness condition $T:A_{j}\overset{c}{%
\rightarrow }B_{j}$ for \textit{both} values $0$ and $1$ of $j$. Not only do
we still not know whether in general $T:\vec{A}\overset{c,b}{\rightarrow }%
\vec{B}$ implies that $T:[A_{0},A_{1}]_{\theta }\overset{c}{\rightarrow }%
[B_{0},B_{1}]_{\theta }$, but we cannot even deduce that $%
T:[A_{0},A_{1}]_{\theta }\overset{c}{\rightarrow }[B_{0},B_{1}]_{\theta }$
when $T$ satisfies the stronger condition $T:\vec{A}\overset{c,c}{%
\rightarrow }\vec{B}$. We conjecture however that $T:\vec{A}\overset{c,c}{%
\rightarrow }\vec{B}$ implies $T:[A_{0},A_{1}]_{\theta }\overset{c}{%
\rightarrow }[B_{0},B_{1}]_{\theta }$ for each $\theta \in (0,1)$ for
arbitrary couples $\vec{A}$ whenever $\vec{B}$ is an arbitrary lattice
couple, i.e., with no requirements of Fatou property or $\sigma $-order
continuity.

\smallskip

\subsection{What if the underlying measure space is not $\sigma $-finite?}

We have assumed throughout this paper that the underlying measure spaces of
our lattice couples are $\sigma $-finite. Could there be an exotic
counterexample in the realm of non $\sigma $-finite measure spaces which
would finally settle the question of whether or not $(*.*)%
\blacktriangleright (*.*)$? This somehow seems unlikely. For a start we can
assert that $\vec{A}\blacktriangleright (*.*)$ also when the underlying
measure space of the lattice couple $\vec{A}$ is not $\sigma $-finite. Let
us note that the proof in \cite{CwK} that $\vec{A}\blacktriangleright (*.*)$
for all lattice couples $\vec{A}$ (Corollary 7 part (c) on p.\ 270) assumes,
albeit quite implicitly, because of the auxiliary results which it uses,
that the underlying measure space is $\sigma $-finite. However it is not
very difficult to obtain appropriate variants of those auxiliary results for
the case of an arbitrary underlying measure space. We plan to provide the
details of that in forthcoming paper(s), where we will also present partial
results showing that $(*.*)\blacktriangleright \vec{B}$ for certain lattice
couples $\vec{B}$ which do not satisfy the hypotheses imposed here,
including some which are defined on non $\sigma $-finite measure spaces.

\smallskip

\section{\label{pps}Appendix-Standard proofs of the standard propositions
\ref{tbsep}, \ref{phir} and \ref{koshy}.}

\textit{Proof of Proposition \ref{tbsep}.} The fact that $(X,d)$ is totally
bounded means that, for each $k\in \Bbb{N}$, there exists a finite set of
points $F_{k}$ such that $X=\bigcup_{y\in F_{k}}B(y,2^{-k})$. The set $%
Y:=\bigcup_{k\in \Bbb{N}}F_{k}$ is of course countable, and, for each $x\in
X $, we have $\inf_{y\in Y}d(x,y)=\inf_{k\in \Bbb{N}}\left( \min_{y\in
F_{k}}d(x,y)\right) =\inf_{k\in \Bbb{N}}2^{-k}=0.$ \qed

\textit{Proof of Proposition \ref{phir}. }Suppose first that, for some $r>0$%
, there exists an infinite set $E$ with the stated properties. If $(X,d)$ is
totally bounded then $X$ is the union of a finite collection of balls each
having radius $r/2$. At least one of these balls must contain infinitely
many elements of $E$. But if two distinct elements $x$ and $y$ of $E$ are in
the same ball of radius $r/2$ then they must satisfy $d(x,y)\le r$. This is
a contradiction which shows that $(X,d)$ is not totally bounded.

For the converse implication, suppose that $(X,d)$ is not totally bounded.
Then there exists some $r>0$ such that $X$ is not the union of any finite
collection of balls of radius $r$. We will construct an infinite sequence $%
\{x_{n}\}_{n\in \Bbb{N}}$ of points in $X$ such that $d(x_{m},x_{n})>r$ for
all $m,n\in \Bbb{N}$ with $m\ne n$. First choose $x_{1}$ to be any point of $%
X$. If $d(x_{1},y)\le r$ for every $y\in X$ then $X=B(x_{1},r)$ which we
know must be false. So there exists some point $x_{2}\in X$ such that $%
d(x_{1},x_{2})>r$. We will now use a recursive procedure to obtain the
points $x_{n}$ for all $n>2$. Suppose we have already obtained $k-1$ points $%
x_{1},x_{2},...,x_{k-1}$ in $X$ such that
\begin{equation}
d(x_{n},x_{m})>r\text{ for all }m,n\in \{1,2,...,k-1\}\text{ with }m\ne n.
\label{grty}
\end{equation}
If it were true that $\min_{n\in \{1,2,...,k-1\}}d(x_{n},y)\le r$ for each $%
y\in X$ then this would imply that $X=\bigcup_{n=1}^{k-1}B(x_{n},r)$, which
we know to be false. Thus there exists at least one point $y_{*}\in X$ which
satisfies
\begin{equation}
\min_{n\in \{1,2,...,k-1\}}d(x_{n},y_{*})>r.  \label{zdde}
\end{equation}
If we choose $x_{k}=y_{*}$ then it follows from (\ref{grty}) and (\ref{zdde}%
) that we now have $d(x_{n},x_{m})>r$ for all $m,n\in \{1,2,...,k\}$ with $%
m\ne n$.

Of course we now take $E$ to be the set of all the points $x_{n}$ and the
proof is complete. \qed

\textit{Proof of Proposition \ref{koshy}. }If $(X,d)$ is not totally bounded
then the infinite sequence $\left\{ x_{n}\right\} _{n\in \Bbb{N}}$ which can
be constructed as in the second part of the proof of Proposition \ref{phir}
clearly cannot have a Cauchy subsequence.

It remains to show the reverse implication. Suppose then that $(X,d)$ is
totally bounded and that $\left\{ x_{n}\right\} _{n\in \Bbb{N}}$ is an
arbitrary sequence in $X$. Let $E$ be the set which contains all elements of
this sequence. If at least one element of $E$ coincides with $x_{n}$ for all
$n$ in some infinite subset $W$ of $\Bbb{N}$, then obviously $%
\{x_{n}\}_{n\in W}$ is a Cauchy subsequence of $\left\{ x_{n}\right\} _{n\in
\Bbb{N}}$ and the proof is complete. Thus we may suppose that $E$ contains
infinitely many elements. We will now construct a sequence of balls $\left\{
B_{k}\right\} _{k\in \Bbb{N}}$ such that the following properties hold for
each $k\in \Bbb{N}$:

(i) $B_{k}$ has radius $2^{-k}$, and

(ii) the set $\bigcap_{j=1}^{k}B_{j}$ contains infinitely many elements of
the set $E$.

To begin this construction we simply use the fact that $X$ is the union of a
finite collection of balls of radius $1$ and so at least one of these balls,
which will be our $B_{1}$, must contain infinitely many elements of $E$. Now
the construction of the $B_{k}$'s for $k>1$ can be done recursively. More
specifically, suppose that we have constructed $\;B_{1},B_{2},.....,B_{m-1}$
with the stated properties (i) and (ii) for $k=1,2,...,m-1$. Again by (a) we
have $X=\bigcup_{p=1}^{N}V_{p}$ for some integer $N$, where each of the sets
$V_{p}$ is a ball of radius $2^{-m}$. Therefore $\bigcap_{j=1}^{m-1}B_{j}=%
\bigcup_{p=1}^{N}\left( V_{p}\cap \bigcap_{j=1}^{m-1}B_{j}\right) $. So, for
at least one value of $p$, the set $V_{p}\cap \bigcap_{j=1}^{m-1}B_{j}$ must
contain infinitely many elements of $E$. We choose $B_{m}$ to be $V_{p}$ for
that value of $p$, and thus we have completed the inductive step which is
required to carry out the construction of $B_{k}$ for all $k\in \Bbb{N}$.

\smallskip Now we can very easily obtain, as required, a Cauchy subsequence $%
\{x_{n_{k}}\}_{k\in \Bbb{N}}$ of the given sequence $\{x_{n}\}_{n\in \Bbb{N}%
} $.\ More explicitly, we are going to construct a strictly increasing
sequence $\left\{ n_{k}\right\} _{k\in \Bbb{N}}$ of positive integers such
that $x_{n_{k}}\in \bigcap_{j=1}^{k}B_{j}$ for each $k\in \Bbb{N}$. First we
choose $n_{1}$ so that $x_{n_{1}}$ is some element of $E$ which is contained
in $B_{1}$. Then we proceed inductively. Suppose that we have already
obtained $n_{1},n_{2},....n_{k-1}$ such that $x_{n_{m}}\in
\bigcap_{j=1}^{m}B_{j}$ for $m=1,2,...,k-1$ and $n_{1}<n_{2}<...<n_{k-1}$.
Since $\bigcap_{j=1}^{k}B_{j}$ contains infinitely many elements of $E$,
there must exist some integer $n_{k}$ which satisfies $n_{k}>n_{k-1}$ and $%
x_{n_{k}}\in $ $\bigcap_{j=1}^{k}B_{j}$. Finally it is easy to check that
the elements $x_{n_{k}}$ which we have constructed in this way for all $k\in
\Bbb{N}$ do indeed form a Cauchy sequence: Given any $\epsilon >0$, choose $%
N $ such that $2^{-N+1}<\epsilon $. Then if $N\le p<q$ we have that $%
x_{n_{p}}$ and $x_{n_{q}}$ are both contained in $\bigcap_{j=1}^{N}B_{j}$
and therefore in $B_{N}$. Since $B_{N}$ is a ball of radius $2^{-N}$ it
follows that $d(x_{n_{p}},x_{n_{q}})<2\cdot 2^{-N}<\epsilon $. \qed

\end{document}